\documentclass{amsart}
\usepackage{amssymb}
\usepackage{amsthm}
\newtheorem{theorem}{Theorem}[section]
\newtheorem{lemma}[theorem]{Lemma}

\theoremstyle{definition}

\theoremstyle{remark}

\begin{document}

\title{Character degree sums in finite nonsolvable groups}

\author{Kay Magaard}
\email{K.Magaard@bham.ac.uk (K. Magaard)}

\author{Hung P. Tong-Viet}
\email{tongviet@maths.bham.ac.uk (H.P. Tong-Viet)}
\address{School of Mathematics,
         University of Birmingham,
        Birmingham, B15 2TT, UK}


\date{\today}


\begin{abstract}
Let $N$ be a minimal normal nonabelian subgroup of a finite group $G.$ We will show that there exists a nontrivial irreducible character of $N$ of degree at least $5$ which is extendible to $G.$ This result will be used to settle two open questions raised by Berkovich and Mann, and Berkovich and Zhmud$\acute{}.$
\end{abstract}

\maketitle


\section{Introduction and notations}
\label{intro}

\noindent
All groups are finite. Let $G$ be a group. Denote by $Irr(G)$ the set of all complex irreducible characters of $G.$ Let $N$ be a normal subgroup of $G.$ Let $\theta\in Irr(N)$ be an irreducible character of $N.$ We say that $\theta$ is extendible to $G$ if there exists $\chi\in Irr(G)$ such that the restriction of $\chi$ to $N$ is $\theta,$ that is $\chi_N=\theta.$ There are many papers devoted to finding a sufficient conditions for $\theta$ to be extendible to $G$ (see Gallagher \cite{Gall}, Gagola \cite{Gago} and \cite[Chapter $8$ and $11$]{Isaacs}). In this paper, we are interested in the existence problem, that is, assume $N$ is a normal subgroup of $G,$ is there any non-trivial irreducible character of $N$ that extends to $G?$
We are mostly concerned with nonsolvable groups. Suppose that  $N$ is a minimal normal nonabelian subgroup of a group $G.$ In \cite[Lemma $5$]{Bia},  it is shown that there exists a nontrivial irreducible character $\theta$ of $N$ which is extendible to $G.$ In Theorem \ref{th1} below, we will show that $\theta$ can be chosen with $\theta(1)\geq 5.$ Using this result, we answer two open problems raised by Berkovich and Mann, and Berkovich and Zhmud$ \acute{}.$
\begin{theorem}\label{th1} Suppose that $N$ is a minimal normal nonabelian subgroup of a group $G.$ Then there exists an irreducible character $\theta$ of $N$ such that $\theta$ is extendible to $G$ with $\theta(1)\geq 5.$
\end{theorem}

Let $T(G)$ be the sum of degrees of complex irreducible characters of $G,$ i.e $T(G)=\sum_{\chi\in Irr(G)}\chi(1).$ Let $k(G)$ be the number of conjugacy classes of $G$ and let $b(G)$ be the largest irreducible character degree of $G.$
Let $N$ be a normal subgroup of $G.$ Denote by $Irr(G,N)$ the set of all complex irreducible characters $\chi$ of $G$ such that $N\not\leq Ker\chi$ and by $T(G,N)$ the corresponding sum of degrees of all characters in $Irr(G,N).$ It is obvious that $Irr(G)=Irr(G/N)\cup Irr(G,N)$ and $T(G)=T(G/N)+T(G,N).$
In \cite[Theorem $8$]{BM}, Y. Berkovich and A. Mann  showed that if $G$ is nonsolvable then $T(G)>2|G:G'|$ and they asked whether or not $T(G)>2T(G/N),$ where $N$ is a nonsolvable normal subgroup of $G.$ Here is our first result.

\begin{theorem}\label{th2} Let $N$ be a nonsolvable normal subgroup of a group $G.$ Then $T(G)\geq 6T(G/N).$
\end{theorem}
This settles Question $4$ in \cite{BM} or Problem $138$ in \cite{BZ}. By Schwarz inequality, it is easy to see that $T(G)^2\leq |G|k(G).$ Hence $T(G)$ can be used to estimate $k(G).$ Moreover, as $T(G)\leq k(G)b(G),$ we can also get a lower bound for $b(G)$ in terms of $T(G)$ and $k(G).$ The reason that we are interested in the character degree sums comes from a question of Jan Saxl (\cite[Problem $9.56$]{MK}) which asked for the classification of groups in which the square of every irreducible character is multiplicity free.
In fact if we could prove that $T(G/F(G))<b(G/F(G))^2,$ provided $G$ is nonsolvable, where $F(G)$ is the Fitting subgroup of $G,$  then such a group in Jan Saxl's question is solvable. This will limit the possibilities for such groups. The proof for this fact is quite straight forward. Let $G$ be a minimal counter-example to the assertion that the square of every irreducible character of $G$ is multiplicity free but $G$ is not solvable. We first observe that if $N\unlhd G$ and $\chi\in Irr(G/N)$ then $\chi\in Irr(G)$ and every irreducible constituent of $\chi^2$ in $G$ is also an irreducible character of $G/N$ so that $\chi^2$ is  multiplicity free in $G/N$ since it is multiplicity free in $G.$ Thus $G/N$ satisfies Saxl's condition. Secondly, for any $\chi\in Irr(G),$ as $\chi^2$ is multiplicity free, it follows that $\chi(1)^2=\chi^2(1)\leq T(G).$ Now combining these two observations for the quotient group $G/F(G),$ we obtain $b(G/F(G))^2\leq T(G/F(G)),$ where $b(G/F(G))$ is the largest character degree of $G/F(G).$ As $G$ is nonsolvable and the Fitting subgroup $F(G)$ is solvable, $G/F(G)$ is nonsolvable. Then the inequality mentioned above would provide a contradiction.

Denote by $T_1(G)$ the sum of degrees of nonlinear irreducible characters of $G.$ Let $Irr_2(G)=\{\chi\in Irr(G)\:|\:\chi(1)>2\}$ and let $T_2(G)$ be the sum of degrees of characters in $Irr_2(G).$ Observe that if $G$ does not have any irreducible characters of degree $2$ then $T_1(G)=T_2(G),$ for example, this is the case if $G$ is a nonabelian simple group.  The following result is a generalization of \cite[Theorem $8$]{BM}.

\begin{theorem}\label{th3} If $G$ is nonsolvable then $T_2(G)\geq 5|G:G'|.$
\end{theorem}

It is well known that a group $G$ is abelian if and only if $T(G)=k(G).$ The following theorem shows that the structure of $G$ is very restricted when $T(G)$ is small in terms of  $k(G).$
\begin{theorem}\label{th4} If $T(G)\leq 2k(G)$ then $G$ is solvable.
\end{theorem}
This gives a positive answer to Problem $24$ in \cite{BZ}. We note that this property does not characterize the solvability of groups. In fact, let $G\cong 3^2:2S_4,$ which is a maximal parabolic subgroup of $PSL(3,3).$ We have  $T(G)=50,$ $k(G)=11,$ $T(G)>4k(G)$ and $G$ is solvable. Now if $G\cong A_5,$ then $T(G)=16, k(G)=5,$ $T(G)<4k(G)$ and $G$ is nonsolvable. We conjecture that a group $G$ is solvable provided that $T(G)\leq 3k(G).$

\section{Preliminaries}

\begin{lemma}\label{lem1} Let $T$ be a non-abelian simple group. Then there exists a nontrivial irreducible character $\varphi$ of $T$ that extends to $Aut(T)$ with $\varphi(1)\geq 5.$
\end{lemma}
This is essentially Lemma $4.2$ in \cite{Moreto} or \cite[Theorems $2,3,4$]{Bia}. However, the fact that $\varphi(1)\geq 5$ is not explicitly stated there so that we will give a proof for completeness.

\begin{proof}
According to the Theorem of Classification of Finite Simple Groups, every nonabelian simple group is isomorphic to the alternating group of degree $n\geq 7,$ a sporadic group or a finite group of Lie type. We will consider the Tits group ${}^2F_4(2)$ as a sporadic rather than a finite group of Lie type. For alternating group $A_n, n\geq 7,$ the irreducible character $\varphi$ corresponding to the partition $(n-1,1)$ extends to $S_n$ and $\varphi(1)=n-1\geq 6.$ If $T$ is a sporadic group, Tits group or $A_5,$ by inspecting \cite{ATLAS}, we can see that there exists an irreducible character $\varphi$ of $T$ that extends to $Aut(T)$ with $\varphi(1)\geq 5.$ Finally assume $T$ is a finite group of Lie type defined over a field of size $q=p^f,$ where $p$ is prime. Choose $\varphi$ to be the Steinberg character of $T$ of degree $|T|_p,$ the order of the $p$-Sylow subgroup of $T.$ Then $\varphi$ is extendible to $Aut(T)$ (see \cite{Feit}). Moreover, we can easily check that $|T|_p> 5$ provided that $T\not\cong L_2(4)\cong L_2(5)\cong A_5.$ Thus $\varphi(1)\geq 5.$    This completes the proof.
\end{proof}

\begin{lemma}\emph{(Gallagher \cite[Corollary $6.17$]{Isaacs}).}\label{lem2} Let $N\unlhd G$ and let $\chi\in Irr(G)$ be such that $\chi_N=\vartheta\in Irr(N).$ Then the characters $\beta \chi$ for $\beta\in Irr(G/N)$ are irreducible, distinct for distinct $\beta$ and are all of the irreducible constituents of $\vartheta^G.$
\end{lemma}



\section{Proof of the Main Results}
\noindent
{\bf Proof of Theorem \ref{th1}.}
 Since $N$ is a minimal normal nonabelian subgroup of $G,$ there exists a nonabelian simple group $T$ such that $N=T_1\times T_2\times\cdots\times T_k,$ where $T_i\cong T, i=1,\cdots, k.$ Let $\varphi$ be an irreducible character of $T$ obtained from Lemma \ref{lem1} and let $\theta=\varphi\times \varphi\times\cdots\times\varphi.$ By \cite[Lemma $5$]{Bia}, $\theta\in Irr(N)$ and it is extendible to $G.$ As $\varphi(1)\geq 5,$ we have  $\theta(1)=\varphi(1)^k\geq 5.$  The proof is now complete. $\hfill\square$\\

\noindent
{\bf Proof of Theorem \ref{th2}.} We argue by induction on the order of $G.$
Assume first that $N$ is a minimal normal subgroup of $G.$ By Theorem \ref{th1}, there exists an irreducible character $\varphi$ of $N$ which extends to an irreducible character $\chi$ of $G$ with $\chi(1)\geq 5.$  Now by Lemma \ref{lem2}, there is an injective map from $Irr(G/N)$ to $Irr(G,N)$ which maps $\beta\in Irr(G/N)$ to $\beta\chi\in Irr(G,N),$ so that $\chi(1)T(G/N)\leq T(G,N).$ Therefore $$T(G)=T(G/N)+T(G,N)\geq (1+\chi(1))T(G/N)\geq 6 T(G/N).$$

Let $K$ be a minimal normal subgroup of $G$ which is contained in $N.$ Then $K$ is a proper subgroup of $N.$
If $K$ is solvable, then $N/K$ is nonsolvable, and by inductive hypothesis, we have $$T(G)\geq T(G/K)\geq 6T((G/K)/(N/K))=6T(G/N).$$
If $K$ is nonsolvable, then we can apply the result proved in the first paragraph to deduce that $T(G)\geq 6T(G/K).$ As $K\leq N,$ we have $T(G/K)\geq T(G/N)$ so that $T(G)\geq 6T(G/K)\geq 6T(G/N).$
The proof is now complete. $ \hfill\square$\\

\noindent
{\bf Proof of Theorem \ref{th3}.}
Let $N$ be the last term of the derived series of $G$ and let $K$ be maximal among the normal subgroups of $G$ that are contained in $N.$ Since $N=N'\leq G',$ it suffices to prove the result for $G/K$ so that we can assume that $K=1$ and hence $N$ is a minimal normal nonabelian subgroup of $G.$ By Theorem \ref{th1}, there exists an irreducible character $\chi\in Irr(G),$ with $\chi(1)\geq 5,$ and
$\chi_N=\varphi\in Irr(N).$ Let $\psi=\chi_{G'}.$ As $\psi_N=\varphi\in Irr(N),$ it follows that $\psi\in Irr(G')$ and hence $\chi_{G'}=\psi\in Irr(G')$ with $\chi\in Irr(G)$ and $\chi(1)\geq 5.$ Now by Lemma \ref{lem2}, there is an injective map from $Irr(G/G')$ to ${Irr}_2(G)$
which maps $\beta\in Irr(G/G')$ to $\beta\chi\in {Irr}_2(G),$ so that $\chi(1)|G:G'|\leq T_2(G).$ Thus $T_2(G)\geq 5|G:G'|.$ This finishes the proof. $ \hfill\square$ \\

\noindent
{\bf Proof of Theorem \ref{th4}.}
 By way of contradiction, assume that $G$ is nonsolvable. Let $a$ be the number of linear characters of $G,$ let $b$ be the number of irreducible characters of $G$ of degree $2$ and finally let $c$ be the number of irreducible characters of degree greater than $2.$ We have
\begin{equation}\label{eqn1}
    k(G)=a+b+c
\end{equation}

\begin{equation}\label{eqn2}
   T(G)=a+2b+T_2(G)
\end{equation}

\begin{equation}\label{eqn3}
   T(G)\geq a+2b+3c
\end{equation}
Since $T(G)\leq 2k(G),$ it follows from  (\ref{eqn1}) and (\ref{eqn3}) that $$a+2b+3c\leq 2a+2b+2c$$ and hence
\begin{equation}\label{eqn4}
    c\leq a
\end{equation}
Since $T(G)\leq 2k(G),$ it follows from (\ref{eqn1}) and (\ref{eqn2}) that $$a+2b+T_2(G)\leq 2a+2b+2c$$ and so
$$T_2(G)\leq a+2c.$$ Combining with (\ref{eqn4}), we obtain
\begin{equation}\label{eqn5} T_2(G)\leq 3a=3|G:G'|.
\end{equation}
However, this contradicts  Theorem \ref{th3}. Thus $G$ must be solvable. This completes the proof.  $ \hfill\square$

\subsection*{Acknowledgment} The authors would like to thank the referee
for the helpful comments and suggestions. These have improved our
exposition and have shortened the proof of Theorem \ref{th1}. The second author
is supported by the Leverhulme Trust and is grateful to Dr. Paul Flavell
for his help and support.

\bibliographystyle{amsplain}

\end{document}